\documentclass[final, 1p]{elsarticle}

\usepackage{amsmath}
\usepackage{amsfonts}
\usepackage{amssymb}
\usepackage{mathrsfs}
\usepackage{fancybox}
\usepackage{makeidx}
\usepackage{graphicx}
\usepackage{subfigure}
\usepackage{epsfig}
\usepackage{verbatim}
\usepackage{hyperref}
\usepackage{cases}
\usepackage{booktabs}
\usepackage{upgreek}
\usepackage{tabularx,multirow, color}
\usepackage[active]{srcltx}
\usepackage{citesort}
\usepackage{ulem}   
\let \oldmarginpar \marginpar
\renewcommand{\marginpar}[1]{\oldmarginpar{\color{red}{#1}}}

\newcounter{algleo}
\newlength{\lefttab}
\newlength{\numberoffset}
\newenvironment{algleo}%
  {\trivlist
   \topsep=0pt\parsep=0pt\itemsep=0pt
   \def\li{\item\refstepcounter{algleo}\makebox[0.8em][r]{\thealgleo\hspace{\numberoffset}}
       \hangafter1\hangindent1.8em\noindent}%
   \def\linonumber{\item\makebox[0.8em][r]{\hspace{\numberoffset}}
       \hangafter1\hangindent1.8em\noindent}%
   \addtolength{\lefttab}{1.25em}
   \addtolength{\numberoffset}{1.25em}
   \leftskip=\lefttab}%
  {\endtrivlist}

\begin{document}

\begin{frontmatter}

\title{A SVD accelerated kernel-independent fast multipole method and its application to BEM}

\author[]{Yanchuang Cao}
\ead{caoyanch@126.com}

\author[]{Lihua Wen\corref{cor1}}
\ead{lhwen@nwpu.edu.cn}

\author[]{Junjie Rong}
\ead{rxrjj@126.com}

\cortext[cor1]{Corresponding author}
\address{College of Astronautics, Northwestern Polytechnical University, Xi'an 710072, P. R. China}

\begin{abstract}
The kernel-independent fast multipole method (KIFMM) proposed in \citep{kifmm} is of almost linear complexity. In the original KIFMM the time-consuming M2L translations are accelerated by FFT. However, when more equivalent points are used to achieve higher accuracy, the efficiency of the FFT approach tends to be lower because more auxiliary volume grid points have to be added. In this paper, all the translations of the KIFMM are accelerated by using the singular value decomposition (SVD) based on the low-rank property of the translating matrices. The acceleration of M2L is realized by first transforming the associated translating matrices into more compact form, and then using low-rank approximations. By using the transform matrices for M2L, the orders of the translating matrices in upward and downward passes are also reduced. The improved KIFMM is then applied to accelerate BEM. The performance of the proposed algorithms are demonstrated by three examples. Numerical results show that, compared with the original KIFMM, the present method can reduce about 40\% of the iterating time and 25\% of the memory requirement.
\end{abstract}


\end{frontmatter}

\section{Introduction}

The boundary element method (BEM) has become a promising numerical method in computational science and engineering. Despite many unique advantages, like the dimension reduction, high accuracy and suitability for treating infinite domain problems, a major disadvantage of the BEM is its dense system matrix which solution cost is prohibitive in large-scale problems.
During the past three decades, several acceleration methods
have been proposed to circumvent this disadvantage. Representative examples are the fast multipole method (FMM)\cite{fmm},
wavelet compression method\cite{wbem}, $\mathcal{H}$-matrix\cite{Hmatrix},
adaptive cross approximation (ACA)\cite{ACA}, pre-corrected FFT
\cite{pFFT}, etc. Among them the FMM is no doubt the most outstanding one.

The conventional FMM is originally proposed to accelerate the $N$-body simulations, which requires the analytical expansions of
the kernel functions. This poses a severe limitation on its applications to many problems where the analytical expansions are hard to be obtained. Besides, the kernel-tailored expansion makes it difficult to develop a universal FMM
code for real-world applications. To overcome this drawback, the kernel-independent FMM (KIFMM) has been proposed in the past decade
\cite{kifmm, bbFMM, gfmm}. A salient feature of the KIFMM is that the expansion of the kernel function is no longer required.
Instead, only a user-defined function for kernel value evaluation is needed; the structure of the FMM acceleration algorithm is in common for many typical problems.

In this paper, the KIFMM proposed by Ying et al \cite{kifmm} is concerned. This method uses equivalent densities
in lieu of the analytical expansions. It provides a unified framework for fast summations with the Laplace, Stokes, Navier and similar kernel functions. Due to its ease-of-use  and high efficiency,
it has attracted the attention of many researchers \cite{kifmmapp1, kifmmapp2, kifmmapp3}.

The moment-to-local (M2L) translation is the most time-consuming part of the FMM \cite{pwFMM,
fftFMM, gfmm, ecFMM, bbFMM, expFMM}. In the KIFMM \cite{kifmm} the M2L translation is accelerated by the fast Fourier transform (FFT), leading to $\mathcal{O}(p^3\log p)$ computational complexity, where $p$ is the number of equivalent points along the cube side. However, one should note that the efficiency of the FFT approach tends to become lower when $p$ increases. This is because the equivalent points lie only on the boundary of each box, while to use the FFT Cartesian grid points interior the box must be considered as well. In this paper, the M2L in KIFMM is compressed and accelerated using the singular value decomposition (SVD); see Section \ref{improve}. This method is built on the fact that the M2L matrices are typically of very low numerical ranks. Our numerical experiments, including those in Section \ref{numxmp}, show that the proposed method is more efficient than the FFT approach. Another advantage of the SVD accelerating approach is that it is more flexible than the FFT approach, because the later requires the equivalent and check points to be equally spaced while this is not needed in the SVD approach. Moreover, as a byproduct, the orders of the translating matrices in the upward and downward passes can also be reduced by using the transform matrices of M2L, leading to further reduction of the CPU time and memory usage for the upward and downward passes.

The original KIFMM in \cite{kifmm} is designed to accelerate the potential
evaluation for particle simulations. Recently, this original method
was applied to solve boundary integral equations (BIEs) in, e.g., blood flow, molecular electrostatic problems \cite{redbloodcells,
kifmmnys, kifmmbemnys}. It is noticed that the central idea of all those works is to translate the far-field interactions to a particle summation formulation so that the original KIFMM can be used in a \emph{straightforward} manner. For example, Ref. \cite{redbloodcells} deals with large-scale blood flow problem. The velocity of each red blood cell is divided into two components, namely the velocity of a reference point and the relative velocity reflecting the self deformation of the cell. By doing this, the interactions between red blood cells can be formulated into ``particle summations'' corresponding to the reference points for all the blood cells, and thus the KIFMM can be used. In \cite{kifmmnys}, the Nystr\"om method is used to discretize the BIE in order to obtain the particle summation form.

As is well known, the BEM is an dominate numerical method for solving BIEs and has profound applications in engineer. In this paper, the KIFMM is used to accelerate the BEM. This work is nontrivial since the KIFMM can not be straightly used in BEM due to the presence of elements, let alone to maintain the accuracy and efficiency. For example, the equivalent and check surfaces are crucial components of the KIFMM. In the original KIFMM these surfaces can be set as the surfaces of each cube. However, in BEM setting this choice would deteriorate the accuracy, because the boundary elements belonging to a cube can often extrude from the cube; see Section \ref{S-S-equ_chk_sf} for the details in choosing those surfaces.

\section{Basic idea of the KIFMM} \label{kifmmrev}

The KIFMM was proposed in \cite{kifmm} to solve the potential problems for particles. Here its framework is briefly reviewed.

\subsection{Setup} \label{S-2-setup}

Assume that there are $N$ source densities $\{q_i\}$ located at $N$ points
$\{\mathbf{y}_i\}$. Then the induced field potential
$\{p_i\}$ at points $\{\mathbf{x}_i\}$ is given by
\begin{equation} \label{eq_eval}
	p_i = p(\mathbf{x}_i) = \sum_{j=1}^N G(\mathbf{x}_i, \mathbf{y}_j)
	q(\mathbf{y}_j)
	= \sum_{j=1}^N G_{ij} q_j, \quad i = 1, 2, \cdots, N,
\end{equation}
where, $G(\mathbf{x}, \mathbf{y})$ is the kernel function, which can be
of the single layer, double layer or other layers. The complexity is obviously
$\mathcal{O}(N^2)$ if the potentials are evaluated naively by a order $N$
matrix-vector multiplication. By using the FMMs this complexity can be
reduced to $\mathcal{O}(N)$.

The central to all the FMMs lies in the low-rank approximation (LRA) of the submatrices representing the far-field interactions.
The efficient realization of the LRA relies on a spatial tree structure. To construct the tree,
all the particles are first included into a root level cube. Then the cube is equally
divided into 8 cubes, generating the cubes in the next level.
This subdivision is continued until the particles contained in each leaf cube is no more than a
predetermined number $s$.

For each cube $C$, let $\mathscr{N}^C$ denote its near field which consists of cubes in the same level that share at least one vertex with $C$; the union of the other cubes is defined to be its far field $\mathscr{F}^C$.
Let $B$ denote the parent cube of $C$, then the interaction field of $C$ is defined as
$\mathscr{I}^C = \mathscr{F}^C \backslash \mathscr{F}^B$.
Let $\mathbf{y}^{C, \mathrm{u}}$ denote the \emph{upward
equivalent surface} corresponding to cube $C$, $\mathbf{x}^{C, \mathrm{u}}$ denote the \emph{upward check
surface}, $\mathbf{y}^{C, \mathrm{d}}$ denote the \emph{downward
equivalent surface} and  $\mathbf{x}^{C, \mathrm{d}}$ denote the \emph{downward check
surface} \cite{kifmm}. To guarantee the existence of the equivalent densities and check potentials, these surfaces must satisfy the following conditions:
\begin{enumerate}
	\item $\mathbf{y}^{C, \mathrm{u}}$ and $\mathbf{x}^{C, \mathrm{u}}$
		lie between $C$ and $\mathscr{F}^C$;
		$\mathbf{x}^{C, \mathrm{u}}$ encloses $\mathbf{y}^{C, \mathrm{u}}$;
	\item $\mathbf{y}^{C, \mathrm{d}}$ and $\mathbf{x}^{C, \mathrm{d}}$
		lie between $C$ and $\mathscr{F}^C$;
		$\mathbf{y}^{C, \mathrm{d}}$ encloses $\mathbf{x}^{C, \mathrm{d}}$;
	\item $\mathbf{y}^{C, \mathrm{u}}$ encloses $\mathbf{y}^{B,
		\mathrm{u}}$ where $B$ is $C$'s child;
	\item $\mathbf{y}^{C, \mathrm{u}}$ is disjoint from $\mathbf{y}^{B,
		\mathrm{d}}$ for all $B$ in $\mathscr{F}^B$;
	\item $\mathbf{y}^{C, \mathrm{d}}$ lies inside $\mathbf{y}^{B,
		\mathrm{d}}$ where $B$ is $C$'s parent.
\end{enumerate}

The above conditions can be satisfied by choosing all the related surfaces be the boundaries of cubes that are concentric with the cube. For each cube $C$ with side length $2r$,
$\mathbf{y}^{C, \mathrm{u}}$ and $\mathbf{x}^{C, \mathrm{d}}$ can be chosen as the boundary of the cube with halfwidth
$(1+d)r$, $\mathbf{x}^{C, \mathrm{u}}$ and $\mathbf{y}^{C, \mathrm{d}}$ as
the boundary of the cube with halfwidth $(3-2d)r$, where $0 \le d \le \frac{2}{3}$. Therefore, the distance between the equivalent surface and the check surfaces involved in each translation
is no less than $(2-3d)r$. This relation is used in the original KIFMM \cite{kifmm}, with $d$ being of a small value. In this way the equivalent surface and the check surfaces are well-separated, and high accuracy can be obtained.
However, when the KIFMM is applied to BEM, $d$ has to be set larger, or a large part of the source densities on elements belonging to $C$ may extrude from its upward equivalent surface $\mathbf{y}^{C, \mathrm{u}}$, and the sources belonging to cube $C$ can not be well represented by its equivalent densities. Thus the size of the elements should be considered in defining these surfaces for BEM.
See Section \ref{S-S-equ_chk_sf}.

\subsection{Far field translations}

Generally, in a FMM, the potentials induced by the sources in the near field
are computed directly by \eqref{eq_eval}, which is named as \emph{source-to-target} (S2T) translation.
The potentials induced by the sources in the far field are efficiently evaluated by a series of translations, named as \emph{source-to-multipole} (S2M),
\emph{multipole-to-multipole} (M2M), \emph{multipole-to-local} (M2L), \emph{local-to-local} (L2L) and \emph{local-to-target} (L2T) translations.
The main feature of the KIFMM lies in that the above translations are performed using equivalent densities and check potentials,
while in the conventional FMM the translations are performed using the multipole expansions and local expansions. The algorithm for evaluating the potential contribution of far-field sources in KIFMM is as follows.
\begin{enumerate}
    \item \emph{S2M}. The source densities $\mathbf{q}$ in a leaf cube $B$ are translated into its upward equivalent densities $\mathbf{q}^{B, \text{u}}$; that is,
        \begin{equation} \label{eq-s2m}
        \mathbf{q}^{B, \text{u}} = \mathbf{Sq},
        \end{equation}
        with $\mathbf{S}$ being the translating matrix \cite{kifmm}.
    \item \emph{M2M}. The upward equivalent densities $\mathbf{q}^{B, \text{u}}$ of a cube $B$ are transformed to the upward equivalent densities $\mathbf{q}^{C, \text{u}}$ of its parent $C$,
        \begin{equation} \label{eq-m2m}
        \mathbf{q}^{C, \text{u}} = \mathbf{M} \mathbf{q}^{B, \text{u}},
        \end{equation}
        with $\mathbf{M}$ being the translating matrix.
    \item \emph{M2L}. The upward equivalent densities $\mathbf{q}^{C, \text{u}}$ of cube $C$ are translated to the downward check potentials $\mathbf{p}^{D, \text{d}}$ of cube $D \in \mathscr{I}^C$ in its interaction field
        \begin{equation} \label{eq-m2l}
        \mathbf{p}^{D, \text{d}} = \mathbf{K} \mathbf{q}^{C, \text{u}},
        \end{equation}
        where, the translating matrix $\mathbf{K}$ is computed as
        \begin{equation*}
        K_{ij} = G(\mathbf{x}_i, \mathbf{y}_j),
        \end{equation*}
        with $\mathbf{x}_i$ being the $i$-th downward check point of $D$ and $\mathbf{y}_j$ being the $j$-th upward equivalent point of $C$.
    \item \emph{L2L}. The downward check potentials $\mathbf{p}^{D, \text{d}}$ of cube $D$ are translated to the downward check potentials of its child cube $E$,
        \begin{equation} \label{eq-l2l}
        \mathbf{p}^{E, \text{d}} = \mathbf{L} \mathbf{p}^{D, \text{d}},
        \end{equation}
        with $\mathbf{L}$ being the translating matrix.
    \item \emph{L2T}. The downward check potentials $\mathbf{p}^{E, \text{d}}$ of leaf cube $E$ are translated to the potentials $\mathbf{p}$ on the target points in $E$,
        \begin{equation} \label{eq-l2t}
        \mathbf{p} = \mathbf{T} \mathbf{p}^{E, \text{d}},
        \end{equation}
        with $\mathbf{T}$ being the translating matrix.
\end{enumerate}
Combining equations \eqref{eq-s2m}--\eqref{eq-l2t}, the potential $\mathbf{p}$ on the target points in a leaf cube induced by the source densities $\mathbf{q}$ in another leaf cube in its far field can be computed as
\begin{equation} \label{eq-tran0}
    \mathbf{p = TLKMSq}.
\end{equation}

The M2L translation \eqref{eq-m2l} is the most time-consuming step in the KIFMM. It is accelerated by FFT in the original KIFMM \cite{kifmm}. In its implementation auxiliary points must be added inside the upward equivalent surface and the downward check surface, although one only needs the upward equivalent points and the downward check points on the corresponding surfaces. This makes the FFT approach less efficient when the number of the equivalent points and check points are large because the auxiliary points will account for a large proportion. In the next section, a SVD approach is proposed to accelerate the M2L translations as well as other translations.

\section{SVD-based acceleration for translations} \label{improve}

In this section, a new SVD-based accelerating technique is proposed, which can compress all the transform matrices in KIFMM, thus both the M2L translation and the upward and downward passes are greatly accelerated.

\subsection{Matrix dimension reduction for M2L} \label{sec-cmpm2l}

In the acceleration for M2L, SVD is applied in two stages. In the first stage, the M2L translating matrices are compressed into more compact forms.

Suppose that the kernel function is translational invariant. The union of unique translating matrices over all cubes in each level forms a set of 316 matrices.
To compress these matrices, first collect them into a fat matrix in
which they are aligned in a single row and a thin matrix in which they
are aligned in a single column
\begin{subequations}
\begin{equation} \label{eq-kfat}
	\mathbf{K}_{\text{fat}} = \left[ \mathbf{K}^{(1)} \quad
	\mathbf{K}^{(2)} \quad \dots \quad \mathbf{K}^{(316)} \right],
\end{equation}
\begin{equation}
	\mathbf{K}_{\text{thin}} = \left[ \mathbf{K}^{(1)}; \quad
	\mathbf{K}^{(2)}; \quad \dots; \quad \mathbf{K}^{(316)} \right],
\end{equation}
\end{subequations}
where $\mathbf{K}^{(i)}$ is the $i$-th translating matrix. Perform SVD for these two matrices
\begin{subequations}
\begin{equation} \label{fat}
	\mathbf{K}_\text{fat} = \mathbf{U} \mathbf{\Sigma} \left[
	{\mathbf{V}^{(1)}}^\mathtt{T} \quad {\mathbf{V}^{(2)}}^\mathtt{T}
	\quad \dots \quad {\mathbf{V}^{(316)}}^\mathtt{T} \right],
\end{equation}
\begin{equation} \label{thin}
	\mathbf{K}_\text{thin} = \left[ \mathbf{Q}^{(1)}; \quad
	\mathbf{Q}^{(2)}; \quad \dots;
	\quad \mathbf{Q}^{(316)} \right] \mathbf{\Lambda} \mathbf{R}^\mathtt{T}.
\end{equation}
\end{subequations}
Notice that in our algorithm, the entities of M2L matrices $\mathbf{K}^
{(i)}$ are the evaluations of single-layer kernel function. In most
cases, they are symmetric, ie., ${\mathbf{K}^{(i)}}^\mathtt{T} =
\mathbf{K}^{(i)}$, so \eqref{fat} and \eqref{thin} are just transposes
of each other, and the SVD has to be performed only once.

Consider the translating matrix $\mathbf{K}^{(i)}$ for one translation
\begin{equation}
	\mathbf{U}^\mathtt{T} \mathbf{K}^{(i)} \mathbf{R} = \mathbf{\Sigma}
	{\mathbf{V}^{(i)}}^\mathtt{T} \mathbf{R} = \mathbf{U}^\mathtt{T}
	\mathbf{Q}^{(i)} \mathbf{\Lambda}.
\end{equation}
Obviously, $\mathbf{U}^\mathtt{T} \mathbf{K}^{(i)}
\mathbf{R}$ decays both along the rows and columns as quickly as the singular values in
$\mathbf{\Sigma}$ and $\mathbf{\Lambda}$, thus it can be approximated
by its submatrix $\mathbf{\tilde{U}}^\mathtt{T} \mathbf{K}^{(i)}
\mathbf{\tilde{R}}$, therefore
\begin{equation} \label{eq-cmpM2Ls}
    \mathbf{K = U (U^\mathtt{T} K R) R^\mathtt{T} \approx \tilde{U} (\tilde{U}^\mathtt{T} K \tilde{R}) \tilde{R}^\mathtt{T} = \tilde{U} \tilde{K} \tilde{R}},
\end{equation}
where, $\mathbf{\tilde{U}}$ and $\mathbf{\tilde{R}}$
are the tailored matrices consisted by columns corresponding with dominant
singular values that are not less than $\varepsilon_1 \|\mathbf{K}
_\text{fat}\|_2 = \varepsilon_1 \mathbf{\Sigma}_ {0, 0} = \varepsilon_1
\mathbf{\Lambda}_{0, 0}$, and $\mathbf{\tilde{K}}$ is the compressed translating matrix.
Substituting \eqref{eq-cmpM2Ls} into \eqref{eq-tran0} yeilds
\begin{equation} \label{evalcmp1}
	\mathbf{p} = \mathbf{TL \tilde{U} \tilde{K} {\tilde{R}}^\mathtt{T} MSq}.
\end{equation}
Similar compression scheme was also used in \cite{bbFMM}.

The compression \eqref{eq-cmpM2Ls} is performed for M2L translating matrices at all levels. Let $L$ denote the number of levels, then the computational complexity is $\mathcal{O}(L) \sim \mathcal{O}(\log N)$. It can be reduced for the cases of homogeneous kernels. Assume that kernel function $G(\mathbf{x}, \mathbf{y})$ is homogeneous of degree of $m$, that is, $G(\alpha \mathbf{x}, \alpha \mathbf{y}) = \alpha^m G(\mathbf{x}, \mathbf{y})$ for any nonzero real $\alpha$. Let $\mathbf{\tilde{K}}_0^{(i)} \, (i = 1, 2, \dots, 316)$ be the compressed translating matrices constructed from the interacting cubes that are scaled to have unit halfwidth. Then, the compressed translating matrices on the $l$-th level can be computed efficiently by scaling
\begin{equation}
    \mathbf{\tilde{K}}_l^{(i)} = \left( \frac{r_0}{2^l} \right)^m \mathbf{\tilde{K}}_0^{(i)},
\end{equation}
where, $r_0$ is the halfwidth of the root cube in the octree.
Therefore only $\mathbf{\tilde{K}}_0^{(i)}$ has to be computed in the compressing procedure, and the computational complexity can be reduced into $\mathcal{O}(1)$.

The threshold $\varepsilon_1$ affects the balance between the computational cost and the accuracy of the algorithm. The induced error in each M2L translation is of order $\varepsilon_1$, and the total error is approximately $L \varepsilon_1$ \cite{kifmm}. In order to maintain the error decreasing rate of BEM with piecewise constant element, $L\varepsilon_1$ should decrease by a factor of 2 with each mesh refinement
$$
L \varepsilon_1 \sim 2^{-L}.
$$
In this paper, $\varepsilon_1$ is chosen by
\begin{equation} \label{eq-vareps1}
    \varepsilon_1 = C_1 \frac{2^{-L}}{L},
\end{equation}
where, $C_1$ is a constant coefficient.

\subsection{Further acceleration for M2L} \label{subsec-lra}

\begin{figure}[h]
  \centering
  \includegraphics[width = 0.4\textwidth, angle = 0]{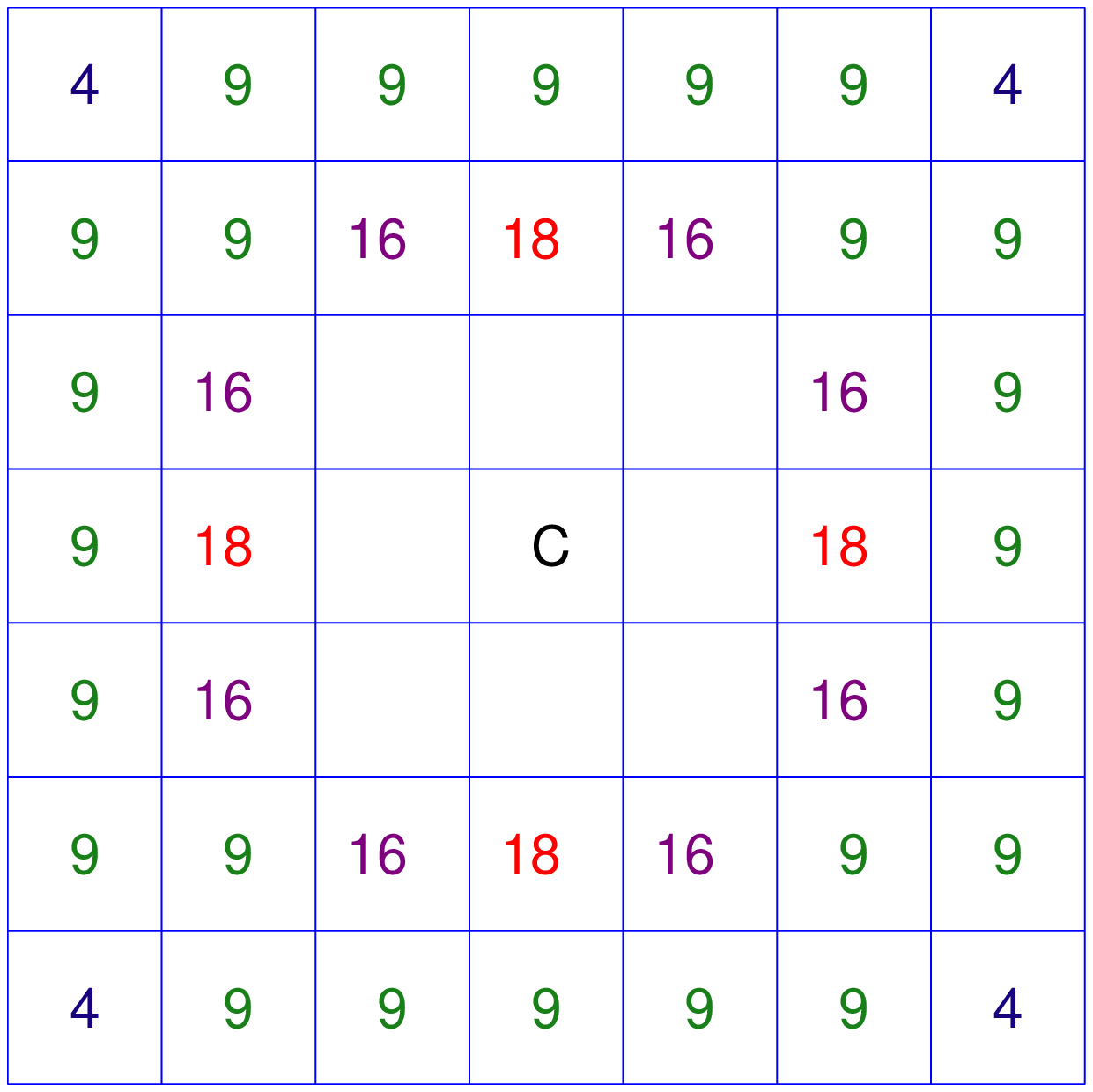}
  \caption{Numerical rank distribution of M2L matrices $\mathbf{\tilde{K}}_{84 \times 84}$ in numerical example \ref{subsec-electrostatic} with $N=2097152, p=8, C_1=0.1, C_2=100$.}
  \label{ranksM2L}
\end{figure}

After the dimension reduction, it is found that most of the compressed M2L matrices $\mathbf{\tilde{K}}$ are still of low numerical ranks. For example, figure \ref{ranksM2L} illustrates the rank distribution of the interaction field of a cube $C$ used in numerical example \ref{subsec-electrostatic} with $N=2097152, p=8, C_1=0.1, C_2=100$.
The dimension of the original translating matrix is 296. After compression using $\mathbf{\tilde{U}}$ and $\mathbf{\tilde{R}}$ the dimension is reduced to 84. However, the figure clearly shows that the actual numerical ranks of the matrices are still much lower than 84.
This fact indicates that the computational cost of M2L can be further reduced by using the low rank decomposition of matrices $\tilde{\mathbf{K}}$. Here the low rank decomposition is computed by SVD, so that optimal rank can be obtained. Since the number of the translating matrices is $\mathcal {O}(1)$, this computational overhead is small.

Consider the low rank approximations of the scaled matrices $\mathbf{\tilde{K}}_0$ for translational invariant and homogeneous kernels.
Compute the SVD for each M2L matrix $\mathbf{\tilde{K}}_0^{(i)}$,
$$
\mathbf{\tilde{K}}_0^{(i)} = \mathbf{U}_0^{(i)} \mathbf{S}_0^{(i)} \left(\mathbf{Q}_0^{(i)}\right)^\mathtt{T}.
$$
Truncate
the singular values smaller than $\varepsilon_2 \|\mathbf{K}_{0, \text{fat}}\|_2$, and discard the corresponding columns in $\mathbf{U}_0^{(i)}$ and $\mathbf{Q}_0^{(i)}$. Then the M2L translation can be approximated by
\begin{equation}\label{eq-lra}
\begin{split}
    \mathbf{\tilde{p}}^{D, \text{d}}_l &= \sum_{C \in \mathscr{I}^D} \left( \frac{r_0}{2^l} \right)^m \mathbf{\hat{U}}_0^{(i)} (\mathbf{\hat{S}}_0^{(i)}
	\mathbf{\hat{Q}}_0^{(i)}) \mathbf{\tilde{q}}^{C, \text{u}}_l \\
	&= \sum_{C \in \mathscr{I}^D} \left( \frac{r_0}{2^l} \right)^m \mathbf{\hat{U}}_0^{(i)} \mathbf{\hat{V}}_0^{(i)} \mathbf{\tilde{q}}^{C, \text{u}}_l,
\end{split}
\end{equation}
where, $\mathbf{\hat{S}}_0^{(i)}$ is the submatrix of $\mathbf{S}_0^{(i)}$ containing the dominant singular values that are no smaller than $\varepsilon_2 \|\mathbf{K}_{0, \text{fat}}\|_2$; $\mathbf{\hat{U}}_0^{(i)}$ and $\mathbf{\hat{Q}}_0^{(i)}$ are the matrices consisted by the corresponding columns of $\mathbf{U}_0^{(i)}$ and $\mathbf{Q}_0^{(i)}$, respectively; and $\mathbf{\hat{V}}_0^{(i)} = \mathbf{\hat{S}}_0^{(i)} \mathbf{\hat{Q}}_0^{(i)}$.

The error of approximation \eqref{eq-lra} is determined by $\varepsilon_2$.
Denote $\mathbf{\hat{K}}_0^{(i)} = \mathbf{\hat{U}}_0^{(i)} \mathbf{\hat{V}}_0^{(i)}$.
From the truncating scheme, there exists
$$
\| \mathbf{\hat{K}}_0^{(i)} - \mathbf{\tilde{K}}_0^{(i)} \|_2 \le \varepsilon_2
\| \mathbf{K}_{0, \text{fat}} \|_2.
$$
For arbitrary $m \times n$ matrix $\mathbf{A}$, one has $\|\mathbf{A}\|
_\text{max} \le \|\mathbf{A}\|_2 \le \sqrt{mn} \|\mathbf{A}\|_\text{max}$. Thus,
$$
\| \mathbf{\hat{K}}_0^{(i)} - \mathbf{\tilde{K}}_0^{(i)} \|_\text{max} \le
\varepsilon_2 \| \mathbf{K}_\text{fat} \|_2.
$$
Let $\mathbf{\hat{K}}_{0, \text{fat}}$ and $\mathbf{\tilde{K}}_{0, \text{fat}}$ denote the fat matrices for $\mathbf{\hat{K}}_0$ and $\mathbf{\tilde{K}}_0$, respectively, which are constructed similarly as $\mathbf{K}_\text{fat}$ in \eqref{eq-kfat}. It is easy to know that
$$
\| \mathbf{\hat{K}}_{0, \text{fat}} - \mathbf{\tilde{K}}_{0, \text{fat}} \|_\text{max} \le
\varepsilon_2 \| \mathbf{K}_{0, \text{fat}} \|_2.
$$
Since the dimension of $\mathbf{\tilde{K}}_{0, \text{fat}}$ is $\tilde{p} \times 316\tilde{p}$, where $\tilde{p}$ is the dimension of $\mathbf{\tilde{K}}_0$, and
$$
\| \mathbf{\hat{K}}_{0, \text{fat}} - \mathbf{\tilde{K}}_{0, \text{fat}} \|_\text{max} \ge \frac{1}{\sqrt{316\tilde{p}^2}} \| \mathbf{\hat{K}}_{0, \text{fat}} - \mathbf{\tilde{K}}_{0, \text{fat}} \|_2,
$$
one has
$$
\| \mathbf{\hat{K}}_{0, \text{fat}} - \mathbf{\tilde{K}}_{0, \text{fat}} \|_2 \le
\varepsilon_2 \sqrt{316\tilde{p}^2} \| \mathbf{K}_{0, \text{fat}} \|_2.
$$
Therefore, the error introduced by the low rank approximation is ensured to be of same order as $\varepsilon_1$ by letting
\begin{equation*}
\varepsilon_2 \sim \frac{\varepsilon_1}{\sqrt{316\tilde{p}^2}} \sim \frac{\varepsilon_1}{\tilde{p}}.
\end{equation*}
In our scheme, it is defined by
\begin{equation} \label{eq-vareps2}
\varepsilon_2 = C_2 \frac{\varepsilon_1}{\tilde{p}},
\end{equation}
where $C_2$ is a constant coefficient.

\subsection{Acceleration for upward and downward passes} \label{subsec-pass}

The transformation matrices $\mathbf{\tilde{U}}$ and $\mathbf{\tilde{R}}$ can also be used to compress the matrices for upward and downward passes. Since the columns of $\mathbf{\tilde{R}}$ are orthonormal, thus $\mathbf{\tilde{R}^\mathtt{T} \tilde{R} = I}$. The potentials in $\mathscr{I}^B$ generated by the upward equivalent densities $\mathbf{q}^{B, \text{u}}$ can be written as follows
\begin{equation} \label{eq-cmpupeden}
\mathbf{K} \mathbf{q}^{B, \text{u}} \approx \mathbf{\tilde{U} \tilde{K} \tilde{R}^\mathtt{T}} \mathbf{q}^{B, \text{u}}
= \mathbf{\tilde{U} \tilde{K} (\tilde{R}^\mathtt{T} \tilde{R}) \tilde{R}^\mathtt{T}} \mathbf{q}^{B, \text{u}}
= \mathbf{\tilde{U} \tilde{K} \tilde{R}^\mathtt{T}} \mathbf{q}_1^{B, \text{u}},
\end{equation}
where $\mathbf{q}_1^{B, \text{u}} = \mathbf{\tilde{R} \tilde{R}^\mathtt{T}} \mathbf{q}^{B, \text{u}}$ is the projection of $\mathbf{q}^{B, \text{u}}$ to the space spanned by the columns of $\mathbf{\tilde{R}}$. This suggests that $\mathbf{q}_1^{B, \text{u}}$ can approximately reproduce the potential field in $\mathscr{I}^C$ excited by $\mathbf{q}^{B, \text{u}}$. In other words, $\mathbf{q}_1^{B, \text{u}}$ can be taken as the new upward equivalent densities for the potential field in $\mathscr{I}^C$.

Now consider $B$'s parent cube $C$ and its interacting field $\mathscr{I}^C$. Since $\mathscr{I}^C$ lies outside $\mathscr{I}^B$, from potential theory we know that $\mathbf{q}_1^{B, \text{u}}$ can also be used to reproduce the potential field in $\mathscr{I}^C$, ie.,
\begin{equation}
\begin{split}
\mathbf{\tilde{U} \tilde{K} \tilde{R}^\mathtt{T} M} \mathbf{q}^{B, \text{u}}
\approx & \mathbf{\tilde{U} \tilde{K} \tilde{R}^\mathtt{T} M} \mathbf{q}_1^{B, \text{u}} \\
= & \mathbf{\tilde{U} \tilde{K} \tilde{R}^\mathtt{T} M \tilde{R} \tilde{R}^\mathtt{T}} \mathbf{q}^{B, \text{u}} \\
= & \mathbf{\tilde{U} \tilde{K} \tilde{R}^\mathtt{T} M \tilde{R} \tilde{R}^\mathtt{T} S} \mathbf{q} \\
= & \mathbf{\tilde{U} \tilde{K} \tilde{M} \tilde{S}} \mathbf{q},
\end{split}
\end{equation}
where,
$$
\mathbf{\tilde{M} = \tilde{R}^\mathtt{T} M \tilde{R}}
$$
is the new translating matrix for M2M;
$$
\mathbf{\tilde{S} = \tilde{R}^\mathtt{T} S}
$$
is the new translating matrix for S2M.

From the symmetry of the algorithm, ie., the upward pass and the downward pass playing the same role in the algorithm, we know that the downward pass can be transformed by $\mathbf{\tilde{U}}$ similarly. Thus,
$$
\mathbf{\tilde{L} = \tilde{U}^\mathtt{T} L \tilde{U}}
$$
is the new translating matrix for L2L;
$$
\mathbf{\tilde{T} = \tilde{T} \tilde{U}}
$$
is the new translating matrix for L2T.

Since both the transformation matrices $\mathbf{\tilde{U}}$ and $\mathbf{\tilde{R}}$ are thin matrices, the new translating matrices $\mathbf{\tilde{S}}$, $\mathbf{\tilde{M}}$, $\mathbf{\tilde{L}}$ and $\mathbf{\tilde{T}}$ are smaller than their original forms, and thus the computational cost of the upward and downward passes can be reduced.

\section{KIFMM for BEM} \label{kifmbem}

In this section, the KIFMM is applied to accelerate the BEM. One should note that in BEM the sources distribute continuously on the boundary elements instead of
on a group of discrete points in the original KIFMM \cite{kifmm}. Therefore, in the BEM the potential evaluation in
S2T and S2M operations must be performed by integration rather than summation.
More importantly, since the continuous sources are represented by nodal basis
functions, the sources has to be grouped based on the supports of nodal functions,
which would make the definition of equivalent and check surfaces different from original KIFMM.

For clarification in explanation, the single-layer BIE for Laplace problem is considered.
Let $\Omega$ be a bounded domain with boundary $\Gamma$. Given a known potential $f(\mathbf{x})$ on the boundary $\Gamma$,
the source density distribution $q(\mathbf{x})$ satisfies
\begin{equation} \label{ibie}
	\int_\Gamma G(\mathbf{x}, \mathbf{y}) q(\mathbf{y}) \mathrm{d}
	\mathbf{y} = f(\mathbf{x}), \quad \mathbf{x} \in \Gamma,
\end{equation}
where, $G(\mathbf{x}, \mathbf{y}) = {1 / (4 \pi |\mathbf{x}- \mathbf{y}|)}$ is the fundamental solution of the Laplace equation. By partitioning the boundary $\Gamma$ into triangular elements and using the piecewise constant basis functions with the
nodal points on element centroids, the collocation BEM leads to a linear system
$$
\mathbf{Aq = b}
$$
with $\mathbf{q}$ consisting of the source densities on each triangles, $\mathbf{b}$ consisting of the known potentials on each collocation points, and
\begin{equation}
	A_{i, j} = \int_{\triangle_j} G(\mathbf{x}_i, \mathbf{y})
	\chi_j(\mathbf{y}) \mathrm{d}\mathbf{y}.
\end{equation}
where, $\mathbf{x}_i$ is the $i$-th collocation point, $\triangle_j$ is
the $j$-th triangle, and $\chi_j(\mathbf{y})$ is the basis function on
$\triangle_j$. When the system solved by iterative methods, the main computational cost is spent on the matrix-vector multiplication (MVM) of which the complexity is $\mathcal{O}(N^2)$. This complexity can be reduced to $\mathcal{O}(N)$ by
the KIFMM.

\subsection{The equivalent and check surfaces}\label{S-S-equ_chk_sf}

As mentioned above, the definition of the equivalent and check surfaces for BEM has to be different with the original KIFMM for particle simulations, because the equivalent surface must enclose all the sources according to the potential theory \cite{kifmm}.

\begin{figure}[h]
	\centering
	\includegraphics[height=0.3\textwidth, angle = 0]{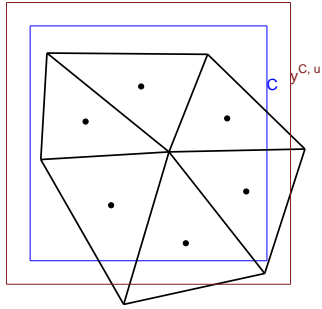}
	\caption{The elements and the upward equivalent surface related to a leaf cube.}
	\label{cube}
\end{figure}

In constructing the octree for BEM, the centroids of elements are used as reference points, and the subdivision process is similar to that in Section \ref{S-2-setup}. Figure \ref{cube} illustrates a leaf cube $C$ in the octree. The union of all the elements whose centroids
lying in $C$ is denoted by $\Gamma(C)$.
One should note that
$\Gamma(C)$ may extrude from $C$. As a result, for BEM the first two
restrictions of surface definition in Section \ref{S-2-setup} have to be
modified as follows:
\begin{enumerate}
	\item $\mathbf{y}^{C, \mathrm{u}}$ and $\mathbf{x}^{C, \mathrm{u}}$
		lie between $\Gamma(C)$
		and $\mathscr{F}^C$; $\mathbf{x}^{C, \mathrm{u}}$ encloses
		$\mathbf{y}^{C, \mathrm{u}}$;
	\item $\mathbf{y}^{C, \mathrm{d}}$ and $\mathbf{x}^{C, \mathrm{d}}$
		lie between $C$ and $\Gamma(\mathscr{F}^C)$, with $\Gamma(\mathscr{F}^C)$ being the union of all elements that belongs to
$\mathscr{F}^C$; $\mathbf{y}^{C, d}$ encloses $\mathbf{x}^{C, \mathrm{d}}$.

\end{enumerate}

The equivalent and check surfaces for BEM are defined similarly with Section \ref{S-2-setup}. However, in order to satisfy the above restrictions, the relative distance between a cube and its upward equivalent surface $d$ has to be chosen large enough so that for each leaf cube, the triangles ``belonging to'' it are enclosed in its upward equivalent surface. For a quasi-uniform element partition, assume that the size of the element is $h$ and each leaf cube contains at most $s$ elements, then the halfwidth of the leaf cubes in the finest level is proportional with $\sqrt{s}h$. The distance between the outmost vertex and the cube surface is no larger than $h$, thus $d$ is of order
\begin{equation*}
d \sim \mathcal {O}\left(\frac{(\sqrt{s}+1)h}{\sqrt{s}h} - 1\right) = \mathcal {O}\left(\frac{1}{\sqrt{s}}\right).
\end{equation*}
So, in this paper $d$ is evaluated as
\begin{equation}
d = C_d \frac{1}{\sqrt{s}},
\end{equation}
where, $C_d$ is user-defined constant. Our numerical experience indicates that $C_d=0.5$ is proper for most problems.

\subsection{S2T and S2M translations}

The equivalent points and check points are sampled on the equivalent and check surfaces, respectively. Then the potentials on the collocation points can be evaluated efficiently by KIFMM, in which the contribution of the near-field sources are evaluated by S2T, and the contribution of the far-field sources are evaluated efficiently using the equivalent densities and check potentials.

The potentials produced by near-field sources are evaluated by S2T translation. In the original KIFMM, it is performed by direct evaluation \eqref{eq_eval}, since the sources distribute on discrete points. However, in BEM the sources distribute on elements, thus these potentials should be evaluated by integration
\begin{equation}
\begin{split}
	p(\mathbf{x}_i) &= \int_{\Gamma(\mathscr{N}^C)} G(\mathbf{x}_i,
		\mathbf{y}) q(\mathbf{y}) \mathrm{d}\mathbf{y} \\
		&= \sum_j q_j \int_{\triangle_j} G(\mathbf{x}_i,
			\mathbf{y}) \chi_j(\mathbf{y})
		\mathrm{d}\mathbf{y}, \quad \triangle_j \in \Gamma(\mathscr{N}^C),
\end{split}
\label{s2teval}
\end{equation}
where $p(\mathbf{x}_i)$ is the check potential on the $i$-th collocation point.

The potentials produced by far-field sources are evaluated by a series of translations, namely S2M, M2M, M2L, L2L and L2T. Among these translations, S2M need to evaluate the upward check potentials produced by the sources belonging to the leaf cube. Similar to S2T translation, this must be implemented by integration as well
\begin{equation}
\begin{split}
	p^{C, \mathrm{u}}(\mathbf{x}) &= \int_{\Gamma(C)} G(\mathbf{x},
		\mathbf{y}) q(\mathbf{y}) \mathrm{d}\mathbf{y} \\
		&= \sum_j q_j \int_{\triangle_j} G(\mathbf{x}, \mathbf{y})
		\chi_j(\mathbf{y}) \mathrm{d}\mathbf{y}, \quad \triangle_j
		\in \Gamma(C),
\end{split}
\label{s2meval}
\end{equation}
where, $p^{C, \mathrm{u}}(\mathbf{x})$ is the upward check potential for leaf cube $C$.

In the above sections, the accelerating algorithm for single layer type boundary integral is introduced. With slight modifications it can be used to accelerate double layer boundary integral. That is, only the integral kernel function in S2T translation and the first step in S2M translation should be replaced into double layer kernel. Therefore, upward equivalent densities should be solved by
\begin{equation} \label{eq-doubles2m}
    \int_{\mathbf{y}^{C, \mathrm{u}}} G(\mathbf{x}, \mathbf{y}) q^{C, \mathrm{u}}(\mathbf{y}) \mathrm{d}\mathbf{y}
    = \int_{\Gamma(C)} \frac{\partial G(\mathbf{x}, \mathbf{y})}{\partial \mathbf{n(y)}} q(\mathbf{y}) \mathrm{d}\mathbf{y}, \quad \text{for any } \mathbf{x} \in \mathbf{x}^{C, \mathrm{u}}.
\end{equation}
Discretized with upward check points and upward equivalent points, a linear system can be achieved, and the upward equivalent densities can be solved. The other steps of the algorithm remains the same with that dealing with single layer boundary integral.

\subsection{The complete algorithm}

In KIFMM for BEM, the discretized sources are grouped into cubes in a octree, then the potentials on collocation points are divided into two parts, namely the contribution of the near-field sources and the contribution of the far-field sources. The former is evaluated by S2T, while the later is evaluated by a series of translations. The complete algorithm for BEM is implemented by the following steps:

\setcounter{algleo}{0}
\begin{algleo}
    \linonumber \textbf{Algorithm} \quad \textsc{SVD accelerated KIFMM for BEM}

    \linonumber \textsc{Step 1 Setup}
    \begin{algleo}
        \li Construct the octree by subdividing the leaf cube recursively.
        \li For each cube $C$, find the cubes in its near field $\mathscr{N}^C$ and interaction field $\mathscr{I}^C$.
        \li Define the equivalent and check surfaces by the method described in Section \ref{S-S-equ_chk_sf}.
        \li Compute and compress the translating matrices by the compressing approach in Section \ref{improve}.
    \end{algleo}

    \linonumber \textsc{Step 2 Upward pass}
    \begin{algleo}
        \li \textbf{for} each leaf cube $C$ in \emph{postorder} traversal of the tree \textbf{do}
        \begin{algleo}
            \li Compute the upward equivalent densities (S2M).
        \end{algleo}
        \li \textbf{end for}
        \li \textbf{for} each non-leaf cube $C$ in \emph{postorder} traversal of the tree \textbf{do}
        \begin{algleo}
            \li Compute the upward equivalent densities (M2M).
        \end{algleo}
        \li \textbf{end for}
    \end{algleo}

    \linonumber \textsc{Step 3 Downward pass}
    \begin{algleo}
        \li \textbf{for} each non-leaf cube $C$ in \emph{preorder} traversal of the tree \textbf{do}
        \begin{algleo}
            \li Add to the downward check potentials produced by the sources in its interaction list (M2L)
            \li Add to the downward check potentials of its child cubes (L2L)
        \end{algleo}
        \li \textbf{end for}
        \li \textbf{for} each leaf cube $C$ in \emph{preorder} traversal of the tree \textbf{do}
        \begin{algleo}
            \li Evaluate the potentials on the collocation points (L2T)
        \end{algleo}
        \li \textbf{end for}
    \end{algleo}

    \linonumber \textsc{Step 4 Near-field interaction}
    \begin{algleo}
        \li \textbf{for} each leaf cube $C$ in \emph{preorder} traversal of the tree \textbf{do}
        \begin{algleo}
            \li Add to the potential the contribution of near field source densities (S2T), which should be evaluated by Eq. \eqref{s2teval}
        \end{algleo}
        \li \textbf{end for}
    \end{algleo}

\end{algleo}

In our method, the definition of the equivalent and check surfaces are different with the original KIFMM. However, this does not affect the computational cost. The total computational complexity of our new KIFMM for BEM remains $\mathcal{O}(N)$.

\section{Numerical Examples} \label{numxmp}

The performance of our SVD-based accelerating technique and the kernel-independent fast multipole BEM
for Laplace BIEs is demonstrated by three numerical examples. The resulting linear systems
are solved by GMRES solver. The algorithms are implemented based on the public \texttt{kifmm3d} code available from \cite{kifmm3dcodes}. All simulations are carried out on a computer
with a Xeon 5440 (3.00 GHz) CPU and 28 GB RAM.

\subsection{Electrostatic problem} \label{subsec-electrostatic}

In this subsection, the electric charge density on an ellipsoidal conductor
is computed by solving Eq. \eqref{ibie}. The ellipsoid can be described by
$(x_1/2)^2 + x_2^2 + (x_3/3)^2 = 1$. The analytic solution can be expressed
analytically using ellipsoidal coordinates.
The convergence tolerance for GMRES solver is
set to be $10^{-6}$. The surface of the ellipsoid is first discretized into $N=512$ triangular elements, then the mesh is refined 6 times. The finest mesh has $N=2097152$ elements.

\begin{table} [h]
	\centering
	\caption{Errors obtained with different $C_1$ and $C_2$}
	\label{table-error}
	\vspace{1em}
	\begin{tabular}{rrrrrrr}
		\hline
		\multirow{2}{*}{$N$}	& \multicolumn{6}{c}{Relative error}	\\
		\cline{2-7}
								& \multirow{2}{*}{FFT}	& $C_1=0.1$		& $C_1=0.5$		& $C_1=0.1$		& $C_1=0.1$		& $C_1=0.1$		\\
								&						& $C_2=0$		& $C_2=0$		& $C_2=10$		& $C_2=100$		& $C_2=500$		\\
		\hline
		\multicolumn{6}{c}{$p=4$}\\
		      512				& 0.069 936				& 0.070 235		& 0.086 715		& 0.071 639		& 0.074 210		& 0.074 210		\\
		    2 048				& 0.032 898				& 0.033 050		& 0.033 227		& 0.033 080		& 0.039 421		& 0.060 588		\\
		    8 192				& 0.014 010				& 0.014 047		& 0.016 236		& 0.014 055		& 0.015 502		& 0.062 169		\\
		   32 768				& 0.006 697				& 0.006 704		& 0.007 845		& 0.006 702		& 0.007 191		& 0.028 167		\\
		  131 072				& 0.003 517				& 0.003 518		& 0.004 236		& 0.003 521		& 0.003 720		& 0.009 063		\\
	      524 288				& 0.002 960				& 0.002 962 	& 0.003 100 	& 0.002 966		& 0.003 098		& 0.004 419		\\
		2 097 152				& 0.004 880				& 0.004 880		& 0.004 887		& 0.004 883		& 0.004 896		& 0.005 296		\\
		\multicolumn{6}{c}{$p=6$}\\
		      512				& 0.069 923				& 0.070 111		& 0.086 805		& 0.071 593		& 0.074 541		& 0.074 541		\\
		    2 048				& 0.032 901				& 0.033 038		& 0.033 270		& 0.033 082		& 0.043 255		& 0.062 104		\\
		    8 192				& 0.014 001				& 0.014 076		& 0.016 217		& 0.014 081		& 0.016 608		& 0.063 492		\\
		   32 768				& 0.006 641				& 0.006 679		& 0.008 849		& 0.006 681		& 0.007 300		& 0.028 298		\\
		  131 072				& 0.003 182				& 0.003 219		& 0.003 983		& 0.003 220		& 0.003 351		& 0.007 635		\\
	      524 288				& 0.001 579				& 0.001 587 	& 0.002 905 	& 0.001 587		& 0.001 607		& 0.002 321		\\
		2 097 152				& 0.000 790				& 0.000 791		& 0.001 227		& 0.000 791		& 0.000 793		& 0.000 932		\\
		\multicolumn{6}{c}{$p=8$}\\
		      512				& 0.069 923				& 0.070 122		& 0.086 831		& 0.071 622		& 0.074 681		& 0.074 681		\\
		    2 048				& 0.032 901				& 0.033 035		& 0.033 299		& 0.033 075		& 0.046 320		& 0.062 675		\\
		    8 192				& 0.014 001				& 0.014 065		& 0.016 208		& 0.014 071		& 0.016 880		& 0.064 226		\\
		   32 768				& 0.006 641				& 0.006 677		& 0.008 673		& 0.006 670		& 0.007 231		& 0.030 980		\\
		  131 072				& 0.003 182				& 0.003 216		& 0.004 141		& 0.003 217		& 0.003 344		& 0.009 128		\\
	      524 288				& 0.001 579				& 0.001 585 	& 0.002 467 	& 0.001 585		& 0.001 602		& 0.002 632		\\
		2 097 152				& 0.000 789				& 0.000 790		& 0.001 014		& 0.000 790		& 0.000 793		& 0.000 988		\\
		\hline
	\end{tabular}
\end{table}

\begin{table} [h]
	\centering
	\caption{CPU times in each iteration $T_\text{it}$ and the total memory usage with different $C_1$ and $C_2$}
	\label{table-cost}
	\vspace{1em}
	\begin{tabular}{r|rrr|rrr}
		\hline
		\multirow{2}{*}{$N$}	& \multicolumn{3}{c|}{$T_\text{it}$ (s)}					& \multicolumn{3}{c}{Memory usage (MB)}	\\
		\cline{2-4}		\cline{5-7}
								& \multirow{2}{*}{FFT}	& $C_1=0.1$		& $C_1=0.1$			& \multirow{2}{*}{FFT}	& $C_1=0.1$		& $C_1=0.1$		 \\
								& 						& $C_2=10$		& $C_2=100$			&						& $C_2=10$		& $C_2=100$		 \\
		\hline
		\multicolumn{7}{c}{$p=4$}\\
		      512				& $\sim$ 0				& $\sim$ 0		& $\sim$ 0			&     2.7				&      1.8		&      1.8		 \\
		    2 048				&  0.02					&  0.01			&  0.01				&     9.3				&      7.1		&      6.7		 \\
		    8 192				&  0.10					&  0.06			&  0.05				&    31.4				&     27.8		&     27.4		 \\
		   32 768				&  0.37					&  0.31			&  0.24				&   121.3				&    116.9		&    116.3		 \\
		  131 072				&  1.57					&  1.50			&  1.10				&   471.1				&    466.4		&    465.3		 \\
	      524 288				&  6.30					&  7.34			&  5.25				& 1 879.6				&  1 892.8		&  1 891.6		 \\
		2 097 152				& 18.25					& 30.12			& 24.32				& 7 504.3				&  7 598.9		&  7 597.7		 \\
		\multicolumn{7}{c}{$p=6$}\\
		      512				&  0.01					& $\sim$ 0		& $\sim$ 0			&      6.1				&      3.0		&      3.0		 \\
		    2 048				&  0.05					&  0.01			&  0.01				&     19.8				&      8.1		&      7.8		 \\
		    8 192				&  0.31					&  0.06			&  0.05				&     55.8				&     28.3		&     27.9		 \\
		   32 768				&  1.11					&  0.30			&  0.23				&    186.4				&    119.5		&    118.9		 \\
		  131 072				&  4.75					&  1.78			&  1.31				&    736.5				&    493.4		&    492.1		 \\
	      524 288				& 13.26					&  8.56			&  6.46				&  2 917.5				&  2 062.5		&  2 060.8		 \\
		2 097 152				& 76.16					& 45.51			& 34.79				& 11 669.0				&  8 861.7		&  8 859.8		 \\
		\multicolumn{7}{c}{$p=8$}\\
		      512				&   0.01				& $\sim$ 0		& $\sim$ 0			&     14.4				&      7.7		&      7.7		 \\
		    2 048				&   0.17				&   0.01		&  0.01				&     44.0				&     12.8		&     12.5		 \\
		    8 192				&   0.70				&   0.06		&  0.04				&    100.8				&     33.0		&     32.6		 \\
		   32 768				&   2.98				&   0.30		&  0.17				&    292.3				&    124.2		&    123.6		 \\
		  131 072				&  12.50				&   1.63		&  1.21				&  1 143.1				&    489.5		&    488.3		 \\
	      524 288				&  49.09				&   8.37		&  4.57				&  4 482.6				&  2 067.1		&  2 065.5		 \\
		2 097 152				& 198.03				&  44.89		& 34.73				& 17 924.3				&  8 866.4		&  8 864.6		 \\
		\hline
	\end{tabular}
\end{table}

The accuracy and efficiency of the present KIFMM BEM are mainly determined by parameters $C_1$ in \eqref{eq-vareps1} and $C_2$ in \eqref{eq-vareps2}.
The translating matrices are independent with the boundary since they are only determined by the position of the equivalent and check points which are defined in the same manner as discussed in section \ref{S-S-equ_chk_sf}. The SVD accelerating approach truncates small singular values of these translating matrices, therefore the induced error by SVD acceleration only depends on $C_1$ and $C_2$ when $p$ is sufficiently large. Consequently, $C_1$ and $C_2$ should keep the same values for various boundary element analyses. From equation \eqref{eq-vareps1} we know that with larger $C_1$, more singular values are discarded, and the translating matrices for M2L and upward and downward passes would be compressed into more compact form, thus the computing time could be reduced lower. While on the other hand the error would become larger. Similar conclusions could be made for $C_2$. Consequently, the choices of $C_1$ and $C_2$ are determined by the tradeoff between the accuracy and the efficiency.

First the influence of $C_1$ on the accuracy of the algorithm is tested. Three cases with $C_1$ being $0$, $0.1$ and $0.5$ are computed. The results corresponding to $C_1=0$ are computed using the original FFT-accelerating scheme in \cite{kifmm}. In all the three cases, $C_2$ is set to be 0.

The resulting errors are listed in the second to fourth columns of Table \ref{table-error}. One can see that when $C_1=0.1$ the errors are almost the same as that computed by FFT-accelerating scheme. However, when $C_1=0.5$ errors for $p=6,\,8$ are increased. This indicates that $C_1=0.1$ is nearly optimal for retaining the accuracy. It is noticed that when $p=4$ the error tends to be larger when the DOF is high. This is because the error of the algorithm is also relevant with $p$ and the depth $L$ of the octree. To get higher accuracy, $p$ has to be increased to reduce the error in each translation; see \cite{kifmm} for the details. The errors with $p = 6$ and $p = 8$ are almost the same, this is because the error is bounded by the discretization precision for the BIE. This also indicates that, for this numerical example, $p = 6$ is sufficient to get the same accuracy with conventional BEM.

The influence of $C_2$ is studied by setting $C_2=10,\,100,\,500$ while $C_1=0.1$. In Table \ref{table-error} it can be seen that for $C_2=10$ and $C_2=100$ the results keep almost the same errors; while for $C_2=500$ the errors increase. Although the errors raise with the increase of $C_2$, the test case indicates that a choice of $C_2$ between 10 and 100 can maintains almost the same accuracy. The CPU times $T_\text{it}$ in each iteration and the total memory usage of the two methods, FFT-accelerating approach and the SVD accelerating approach, are listed in Table \ref{table-cost}. In Table \ref{table-cost}, $T_\text{it}$ can reduce considerably with larger $C_2$, but the memory cost only reduce slightly. The reason is that $C_2$ only controls the accuracy and efficiency of the low-rank approximation for M2L matrices, as discussed in section \ref{subsec-lra}. In this problem, since the kernel is translational invariant and homogeneous, the memory cost for M2L matrices is only of $\mathcal{O}(1)$. Therefore, the memory reduction is negligible.

From Table \ref{table-cost} one can see that the CPU time in the SVD approach can be reduced significantly for large $p$, comparing with the FFT approach, since the CPU time for each iteration in the SVD approach is not sensitive to $p$. For example, the $T_\text{it}$s of the SVD approach for the cases $p=6$ are almost the same as that for $p=8$. This is because when $p$ is large, the size of the compressed translating matrices are mainly determined by the compressing threshold $\varepsilon_1$, and the numerical rank of M2L matrices are only determined by $\varepsilon_2$. Both $\varepsilon_1$ and $\varepsilon_2$ are independent with $p$. However, in the FFT approach more auxiliary points has to be added, which makes the FFT approach less efficient. Besides the CPU time, the memory usage can also be considerably reduced in the SVD approach, since the translating matrices used in S2M and L2T are compressed into more condensed form by the scheme in section \ref{subsec-pass}.

This example shows that the accuracy reduces with the increase of $C_1$ and $C_2$. When $C_1=0.1$ and $C_2=10$ the SVD accelerating approach is much more efficient than FFT-accelerating approach without significantly affecting the accuracy. It is also showed that $p=6$ is sufficient to maintain the accuracy of BEM in this case. With $p=6$, $C_1=0.1$ and $C_2=10$, the CPU time cost in each iteration can be reduced about 40\% and the memory cost can be reduced about 25\% by SVD approach compared with the original FFT approach while maintaining the accuracy of BEM. These parameters will be used in the next numerical examples.

\subsection{A mixed boundary condition problem}

\begin{figure} [h]
	\centering
	\includegraphics[width = 0.6\textwidth, angle = 0]{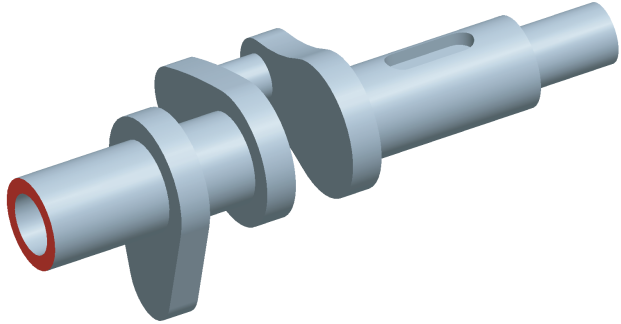}
	\caption{A shaft model with mixed boundary condition.}
	\label{fig-shaft}
\end{figure}

To demonstrate the performance of the SVD accelerating approach for more complicated geometry and boundary condition problems, Laplace equation with mixed boundary conditions on a shaft model illustrated in Figure \ref{fig-shaft} is simulated. The analytical solution is set be to $u = 1/|\mathbf{x_0} - \mathbf{x}|$, with $\mathbf{x_0}$ being outside the computational domain. The potential $u$ is given on the two end surfaces (red surfaces in Figure \ref{fig-shaft}), and the flux $q$ is given on the remainder (gray) surfaces. The converging tolerance for GMRES solver is set to be $10^{-6}$.

The problem is solved by using the KIFMM BEM with the FFT approach and the SVD approach, respectively. The results are reported in Table \ref{table-perf}. It is showed again that the SVD approach can save about 40\% of the iterating time cost and 25\% of the memory cost. The $L^2$-error of $u$ decays linearly with $\mathcal{O}(h)$ and the $L^2$-error of $q$ decays as $\mathcal{O}(\sqrt{h})$. The time consuming in each iteration and the memory consuming increase almost linearly, which indicate that the computational complexity of the method is almost $\mathcal{O}(N)$.

\begin{table} [h]
	\centering
	\caption{Performance of SVD accelerating approach and FFT-accelerating approach}
	\label{table-perf}
	\vspace{1em}
	\begin{tabular}{r|rr|rr|rr|rr}
		\hline
		\multirow{2}{*}{$N$}	& \multicolumn{2}{c|}{Error of $u$}		& \multicolumn{2}{c|}{Error of $q$}		& \multicolumn{2}{c|}{$T_\text{it}$ (s)}		& \multicolumn{2}{c}{Memory (MB)}	\\
		\cline{2-9}
						& FFT			& SVD			& FFT			& SVD			& FFT			& SVD			& FFT			& SVD	\\
		\hline
		   18 048		&  0.021 521	& 0.020 222		& 0.062 718		& 0.066 147		&   0.73		&   0.20		&    115.1		&     77.6	 \\
		   72 518		&  0.012 367	& 0.011 263		& 0.052 828		& 0.053 375		&   2.97		&   0.83		&    396.9		&    281.6	 \\
		  288 768		&  0.004 648	& 0.004 505		& 0.003 173		& 0.003 580		&   7.15		&   4.52		&  1 532.3		&  1 105.4	 \\
		1 156 042		&  0.001 862	& 0.001 866		& 0.002 001		& 0.002 415		&  25.00		&  14.85		&  5 899.9		&  4 533.4	 \\
		\hline
	\end{tabular}
\end{table}

\subsection{Heat conduction problem}

\begin{figure} [ht]
	\centering
	\subfigure[FEM result] {
		\includegraphics[width=0.45\textwidth]{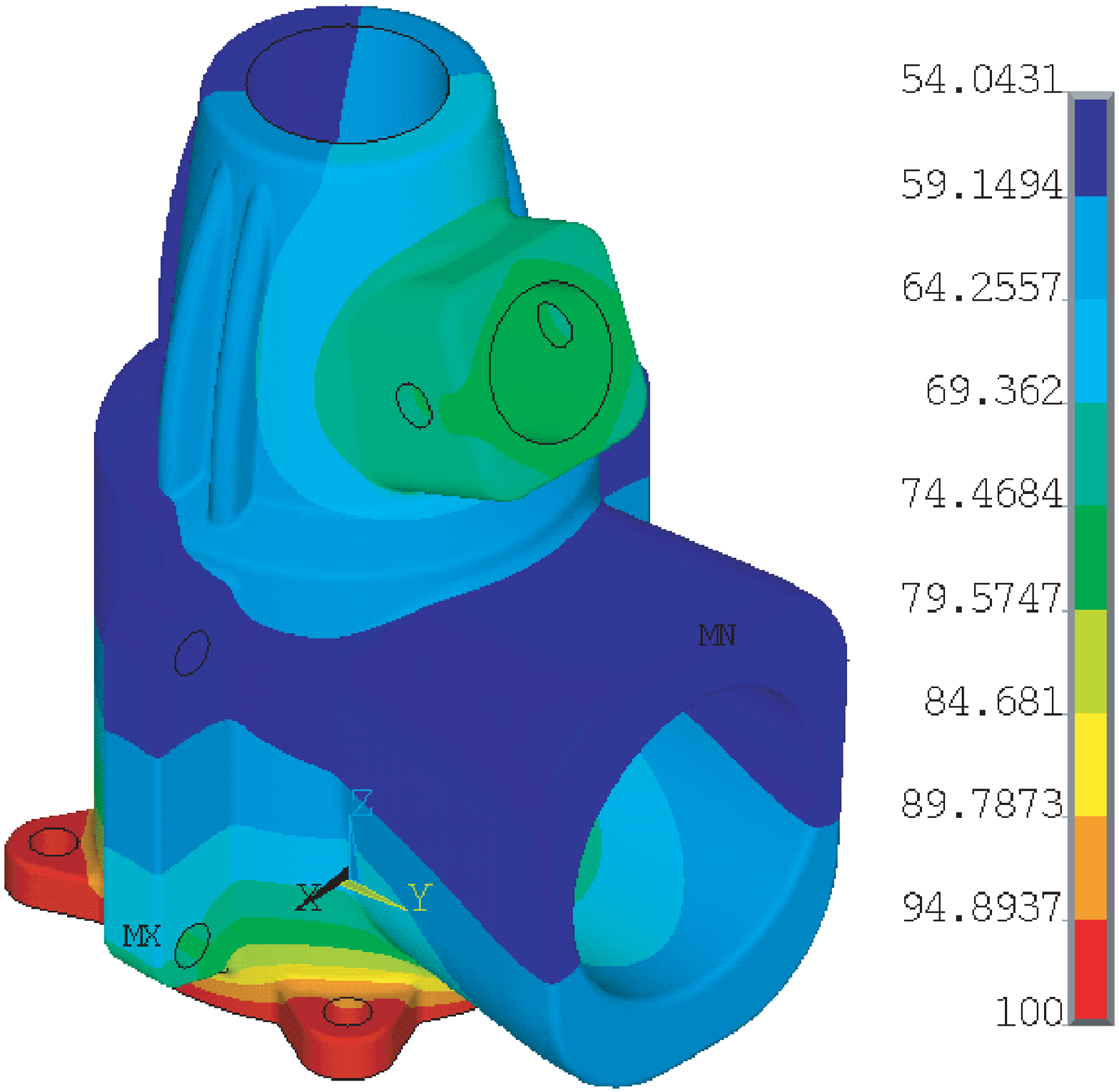}
		\label{engine_fem}
	}
	\subfigure[KIFMM result] {
		\includegraphics[width=0.41\textwidth]{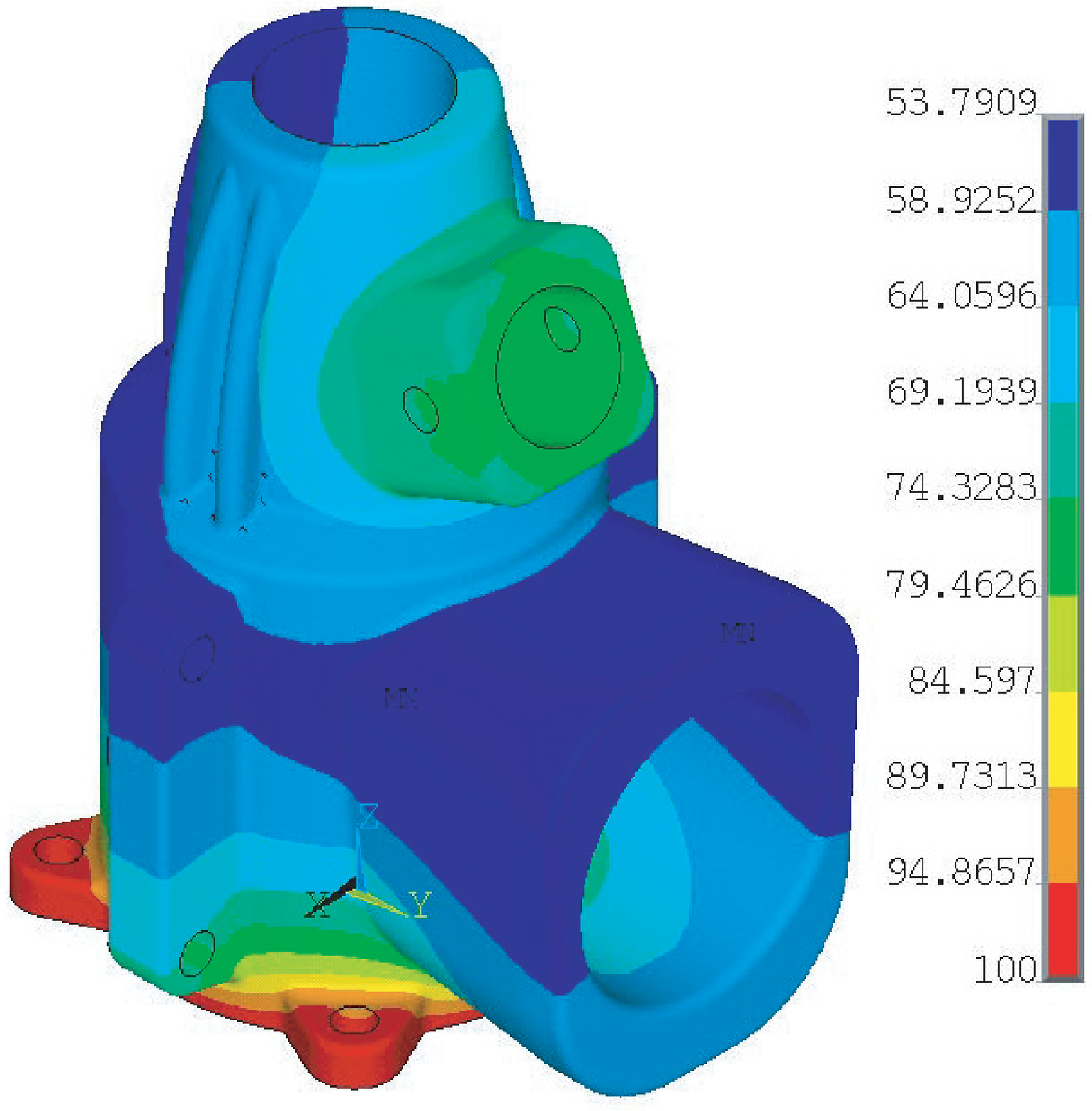}
		\label{engine_fmm}
	}
	\label{engine_result}
	\caption{Temperature distributions of the engine-block model computed by FEM and
		KIFMM BEM.}
\end{figure}

To demonstrate the ability of the present KIFMM BEM for solving real-world problems, a steady-state heat conduction analysis of a engine block is solved here; see figure \ref{engine_result}. The temperature field is governed by the Laplace equation. The conductivity of the engine block is
$\lambda = 80 \mathrm{W/(m \cdot ^\circ C)}$. The temperature of the inner
surface of the oblique tube and the temperature of the bottom surface are
set to be $75 \mathrm{^\circ C}$ and $100 \mathrm{^\circ C}$, respectively.
Convective condition with constant film coefficient $h = 10 \mathrm{W/(m^2 \cdot ^\circ C)}$ and constant
bulk temperature $T_0 = 22^\circ \mathrm{C}$ are applied to the other surfaces. Simulations are performed using a series of meshes with number of elements ranging from 85 680 to nearly 5 million. For comparison, this problem is also solved by finite element method (FEM) with 698317 tetrahedral elements, 1015653 nodes. The converging tolerance for GMRES solver is $10^{-4}$.

The CPU times and memory usage for different
meshes are listed in Table \ref{cost}, where $N_\mathrm{it}$
and $T_\mathrm{it}$ stand for the number of iterations and the CPU time for each iteration, respectively. Again one can see linear behavior of the CPU time and memory requirement. The computed temperature distribution using mesh with 325 774 elements is exhibited in Figure \ref{engine_fmm}. It can be seen that the temperature distribution obtained by the KIFMM BEM agrees very well with that by FEM in figure \ref{engine_fem}.

It is noticed that, with the KIFMM BEM in this paper, the largest model with nearly 5 million DOFs is successfully solved within 5 hours.

\begin{table} [ht]
	\centering
	\caption{CPU times (s) and memory usage (MB) for
		engine-block heat conduction analysis}
	\label{cost}
	\vspace{1em}
	\begin{tabular}{rrrrrr}
		\hline
		$N$           & $N_\mathrm{it}$ 	& $T_\mathrm{total}$	& $T_\mathrm{it}$	& Memory \\
		\hline
		   85 680     &  90		&    582.4		&    3.8		&    502.5\\
		  325 774     &  96		&  2 111.1		&   13.3		&  1 735.5\\
		  900 420     & 100		&  4 251.0		&   27.5		&  4 589.0\\
		1 370 880     & 103		&  5 946.3		&   34.7		&  7 374.8\\
		4 754 670     &  97		& 18 330.4		&  108.2		& 25 021.5\\
		\hline
	\end{tabular}
\end{table}

\section{Conclusion}

The FMM is one of the most successful fast algorithms for BEM acceleration. But it requires the analytical expansion of the kernel function, which makes it difficult to be applied to some complicated problems. Recently, various kernel-independent FMMs were developed to overcome this drawback. Among them the KIFMM proposed in \cite{kifmm} has high efficiency and accuracy, and thus has been extensively used \cite{redbloodcells, kifmmnys, kifmmbemnys}. The KIFMM uses equivalent densities and check potentials instead of analytical expansions to construct the fast algorithm. The time consuming M2L translations are accelerated by using the FFT. However, it is noticed that when more equivalent and check points are sampled to get higher accuracy, the efficiency of the FFT approach tends to be lower because more auxiliary volume grid points have to be added in order to do FFT.

In this paper, the low rank property of the translating matrices in KIFMM is sufficiently exploited by SVD (called SVD approach in this paper) to accelerate all the translations, including the most time-consuming M2L. The acceleration of the M2L translations is carried out in two stages. First the translating matrix is compressed into more compact form, and then it is approximated by low-rank decomposition. By using the compression matrices for M2L, the translating matrices in upward and downward passes can also be compressed into more compact form. Finally, the above improved KIFMM is applied to accelerate BEM, leading to a highly efficient KIFMM BEM for solving large-scale problems.

The accuracy and efficiency of the SVD approach and the KIFMM BEM are demonstrated by three numerical examples. It is shown that, when compared with the FFT-accelerated KIFMM, the SVD approach can reduce about 40\% of the iterating time and 25\% of the total memory requirement. The presented KIFMM BEM is of $\mathcal{O}(N)$ complexity. By using this method Laplace problem with nearly 5 million unknowns can be successfully solved within 5 hours on a Xeon-5440 2.83 GHz CPU and 28 GB RAM.

\section*{Acknowledgements}

This work was supported by the Doctorate Foundation of Northwestern
Polytechnical University under Grant No. CX201220, National Science
Foundations of China under Grants 11074201 and 11102154, and Funds for Doctor Station from the Chinese Ministry of
Education under Grants 20106102120009 and 20116102110006.

\section*{References}
\bibliographystyle{unsrt}
\bibliography{kifmm}

\end{document}